\documentclass[11pt]{article}
\usepackage[table]{xcolor}
\usepackage{graphicx}
\usepackage{amsmath,amsfonts,amssymb}
\usepackage{graphicx}
\usepackage{multirow}
\usepackage{tikz,pgf}
\usepackage{url}
\usepackage{amsmath}
\usepackage{graphicx}
\usepackage{epsfig}
\usepackage{epstopdf}
\usepackage{color}
\usepackage{enumitem}
\def\tto{\;{\lower 1pt \hbox{$\rightarrow$}}\kern -10pt
\hbox{\raise 2pt \hbox{$\rightarrow$}}\;}

\def\Bar{\overline}
\def\ra{\rangle}
\def\la{\langle}
\def\ve{\varepsilon}
\def\epsilon{\varepsilon}
\def\B{\Bbb B}
\def\h{\hfill\Box}
\def\R{\Bbb R}

\def\N{\Bbb N}
\def\ox{\bar{x}}

\def\oy{\bar{y}}

\def\ri{\mbox{\rm ri}\,}

\def\gph{\mbox{\rm gph}\,}
\def\aff{\mbox{\rm aff}\,}
\def\epi{\mbox{\rm epi}\,}

\def\dom{\mbox{\rm dom}\,}
\def\aff{\mbox{\rm aff}\,}

\def\cone{\mbox{\rm cone}}

\def\iri{\mbox{\rm iri}\,}

\def\h{\hfill\square}

\def\emp{\emptyset}

\def\oR{\Bar{\R}}

\def\al{\alpha}

\def\emp{\emptyset}

\def\oR{\Bar{\R}}

\def\al{\alpha}

\def\qri{\mbox{\rm qri}}
\setlist[enumerate,1]{itemsep=0.0ex,parsep=0.5ex,label={\rm(\alph*)},leftmargin=*, align=left}
\newcounter{lk}

\usepackage[colorlinks=true]{hyperref}

\begin{document}
\begin{center}
{\sc\bf Quasi-Relative Interiors for Graphs of Convex Set-Valued Mappings}\\[1ex]
{\sc Dang Van Cuong}\footnote{Department of Mathematics, Faculty of Natural Sciences, Duy Tan University, Da Nang, Vietnam (dvcuong@duytan.edu.vn).},
{\sc Boris S. Mordukhovich}\footnote{Department of Mathematics, Wayne State University, Detroit, Michigan 48202, USA (boris@math.wayne.edu). Research
of this author was partly supported by the USA National Science Foundation under grants DMS-1512846 and
DMS-1808978, by the USA Air Force Office of Scientific Research grant \#15RT04, and by the Australian Research Council under Discovery Project DP-190100555.},
{\sc Nguyen Mau Nam}\footnote{Fariborz Maseeh Department of Mathematics and Statistics, Portland State University, Portland, OR
97207, USA (mnn3@pdx.edu). Research of this author was partly supported by the USA National
Science Foundation under grant DMS-1716057.}\\[2ex]
{\bf Dedicated to Nicolas Hadjisavvas on the occasion of his 65th birthday}
\end{center}
\small{\bf Abstract.} This paper aims at providing further studies of the notion of quasi-relative interior for convex sets. We obtain new formulas for  representing quasi-relative interiors of  convex graphs of set-valued mappings and for convex epigraphs of extended-real-valued functions defined on locally convex topological vector spaces.  We also show that the role, which this notion plays in infinite dimensions and the results obtained in this vein, are similar to those involving relative interior in finite-dimensional spaces. \\[1ex]

\newtheorem{Theorem}{Theorem}[section]
\newtheorem{Proposition}[Theorem]{Proposition}
\newtheorem{Remark}[Theorem]{Remark}
\newtheorem{Lemma}[Theorem]{Lemma}
\newtheorem{Corollary}[Theorem]{Corollary}
\newtheorem{Definition}[Theorem]{Definition}
\newtheorem{Example}[Theorem]{Example}
\renewcommand{\theequation}{\thesection.\arabic{equation}}
\normalsize

\section{Introduction}
\setcounter{equation}{0}

The notion of {\em relative interior} for convex sets in finite-dimensional spaces was largely developed by Rockafellar in his seminal monograph ``Convex Analysis" \cite{r} as a refinement of the classical notion of interior. In contrast to the latter, the relative interior is {\em nonempty} for any nonempty convex subset of $\R^n$. It has been fully recognized by now that results involving relative interior are highly important in many aspects of convex analysis and optimization in finite dimensions. In particular, relative interior is used broadly in developing generalized differentiation theory for convex extended-real-valued functions. For instance, the fundamental result of convex analysis, known as the Moreau-Rockafellar theorem, says that the subdifferential of the sum of two convex functions on $\R^n$ at a point where both functions are finite is represented as the Minkowski sum of the subdifferentials of each functions at references point {\em provided} that intersection of the relative interiors of the domains of these functions is nonempty. Another fundamental result of finite-dimensional geometry tells us the empty intersection of relative interiors of two convex sets is a {\em characterization} of {\em proper separation} of sets. In convex optimization and related topics, the notion of relative interior plays a crucial role in resolving many principal issues of the theory and applications including Fenchel duality, Lagrange duality, optimality conditions, numerical algorithms, etc.; see, e.g., \cite{Bauschke2011,Bertsekas2003,Borwein2000,Boyd2004,g,HP,HU,KK2013,bmn,pr,rw} and the references therein.

It has been recognized for a long time that the classical notion of interior for convex sets and related results in infinite-dimensional spaces, based on convex separation under nonempty interior assumptions, create--besides being rather restrictive in theoretical developments--serious limitations for applications. Among various areas that suffer from such limitations, we mention vector and set optimization as well as general equilibrium theory in economics; see \cite{bao-mor,Flores2,Flores3,ktz,mas,m-book1,m-book} for more discussions and references. In particular, it is well known that the ordering/positive cones in the classical Lebesgue spaces $l^p$ and $L^p$ for $1\le p<\infty$, which naturally appear in economic modeling, have empty interiors, and thus convex separation results conventionally used to establish appropriate price equilibria in basic economic models (e.g., in models of welfare economics) cannot be applied in such frameworks. Moreover, positive cones in the aforementioned and many other infinite-dimensional spaces important for applications have even {\em empty relative interiors}.

All of it provides a strong motivation to seek adequate counterparts of the relative interior notion for convex sets in infinite-dimensions. The major attention in this paper is paid to the notion of {\em quasi-relative interior} for convex sets that was introduced by Borwein and Lewis \cite{bl} and was studied and applied in various publications among which we particularly mention \cite{bao-mor,BG,bot,daniele,durea,ktz,m-book}. As we were informed by Nicolas Hadjisavvas in a private communication, an equivalent notion was introduced independently in his paper with Siegfried Schaible \cite{hs} under the name of {\em inner points}, and this notion can be actually found in Zarantonello \cite{za}.

A key result of \cite{bl,hs} tells us that the quasi-relative interior is {\em nonempty} for every nonempty closed and convex subset of any {\em separable} Banach space. Our main intention here is to study further properties of quasi-relative interior in {\em locally convex topological vector} (LCTV) spaces. We obtain a quasi-relative extension in LCTV spaces of Rockafellar's finite-dimensional theorem on representing the relative interior of convex graphs for set-valued mappings via the relative interiors of the corresponding domain and image sets. The obtained quasi-relative version of such a representation in infinite-dimensional spaces requires a certain {\em quasi-regularity} condition introduced below, which automatically holds in finite dimensions. We also obtain new results on extended relative interiors of convex sets and epigraphs of convex functions on LCTV spaces.

The rest of the paper is organized as follows. Section~2 presents some preliminary material used in the sequel. Section~3 contains the definitions of and discussions on the basic notions of intrinsic relative interior, quasi-relative interior, and quasi-regularity for convex sets in LCTV spaces. Besides other results, we establish here relationships between these properties and the sequential normal compactness of convex sets, which is one the central notions of general infinite-dimensional variational analysis. Section~4 is devoted to the study of quasi-relative interiors for convex graphs of set-valued mappings between LCTV spaces. First we give an alternative proof of Rockafellar's theorem concerning a precise representation of relative interiors for such graphs in finite dimensions. The quasi-interior case of infinite-dimensional spaces is significantly more involved, the results obtained in this section are diverse and require imposing quasi-regularity assumptions. Further developments in this direction for convex sets, set-valued mappings, and extended-real-valued functions in the general LCTV space settings are given in Section~5.

\section{Preliminaries}
\setcounter{equation}{0}

Here we first recall the standard notation and definitions of convex analysis in LCTV spaces; see, e.g., \cite{z}. Then we present, for the reader's convenience, some useful facts on the relative interior of convex sets in finite-dimensional spaces to compare them later with our new results via quasi-relative interiors in general infinite-dimensional settings.

Given a real LCTV space $X$ and its topological dual $X^*$, consider the canonical pairing $\la x^*,x\ra:=x^*(x)$ with $x\in X$ and $x^*\in X^*$. For a nonempty subset $A$ of $X$, define the {\em conic hull} of $A$ by $\cone(A):=\big\{ta\in X\;|\;t\ge 0,\;a\in A\big\}$ and denote the {\em closure} of $A$ by $\Bar{A}$.

Let $\Omega$ be a convex subset of $X$, and let $\ox\in\Omega$. The {\em normal cone} to $\Omega$ at $\ox$ is defined by
\begin{equation}\label{nc}
N(\ox;\Omega):=\big\{x^*\in X\;\big|\;\la x^*,x-\ox\ra\le 0\;\text{ for all }\;x\in\Omega\big\}
\end{equation}
with $N(\ox;\Omega):=\emp$ if $\ox\notin\Omega$. Recall further that two convex sets $\Omega_1,\Omega_2\subset X$ are {\sc properly separated} if there exists $x^*\in X^*$ for which the following two inequality hold:
\begin{equation}\label{ps1}
\sup\big\{\la x^*,w_1\ra\;\big|\;w_1\in\Omega_1\big\}\le\inf\big\{\la x^*,w_2\ra\;\big|\;w_2\in\Omega_2\big\},
\end{equation}
\begin{equation}\label{ps2}
\inf\big\{\la x^*,w_1\ra\;\big|\;w_1\in\Omega_1\big\}<\sup\big\{\la x^*,w_2\ra\;\big|\;w_2\in\Omega_2\big\}.
\end{equation}
Observe that condition \eqref{ps1} can be equivalently rewritten as
\begin{equation*}
\la x^*,w_1\ra\le\la x^*,w_2\ra\;\mbox{\rm whenever }w_1\in\Omega_1,\;w_2\in\Omega_2,
\end{equation*}
while \eqref{ps2} means that there exist $\bar{w}_1\in \Omega_1$ and $\bar{w}_2\in \Omega_2$ such that
\begin{equation*}
\la x^*,\bar{w}_1\ra<\la x^*,\bar{w}_2\ra.
\end{equation*}

Given a nonempty subset $\Omega$ of $X$, the {\em relative interior} of $\Omega$ is defined by
\begin{equation}\label{ri}
\mbox{\rm ri}(\Omega):=\big\{x\in\Omega\;\big|\;\exists\;\mbox{\rm a neighborhood}\, V\, \mbox{\rm of }x\; \mbox{\rm such that }\;V\cap\Bar{\aff}(\Omega)\subset\Omega\big\},
\end{equation}
where $\Bar{\aff}(\Omega)$ denotes the closure of the affine hull of $\Omega$. If $X=\R^n$, the closure operation is not needed in \eqref{ri} since the affine hull $\aff(\Omega)$ is always closed. The following useful result, which is taken from \cite[Theorem~4.7]{bmncal}, provides a proper separation description of the relative interior for convex sets in finite-dimensional spaces.

\begin{Theorem}{\bf(relative interior and proper separation in finite dimensions).}\label{pst1} Let $\Omega_1$ and $\Omega_2$ be two nonempty convex subsets of $\R^n.$ Then $\Omega_1$ and $\Omega_2$ are properly separated if and only if $\ri(\Omega_1)\cap \ri(\Omega_2)=\emptyset.$

\end{Theorem}

Theorem~\ref{pst1} is employed below to provide further characterizations of the relative interior for convex sets. The next theorem combines various known results in this direction (see, e.g., \cite{r}) that are important for understanding the subsequent extensions in infinite dimensions. For completeness, we present a unified and simplified proof of these characterizations.

\begin{Theorem}{\bf(characterizations of relative interior for convex sets in $\R^n$).}\label{pts2} Let $\Omega$ be a nonempty convex set in $\R^n$ and let $\ox\in \R^n$. The following properties are equivalent:\\[1ex]
{\rm(a)} $\ox\in\ri(\Omega)$.\\
{\rm(b)} $\ox\in\Omega$ and for every $x\in\Omega$ with $x\ne\ox$ there exists $u\in\Omega$ such that
$\ox\in(x,u)$, where $(x,u):=\{tx+(1-t)u\in\R^n\;|\;t\in(0,1)\}$ signifies the open interval between the points $x,u$.\\
{\rm(c)} $\ox\in\Omega$ and $\cone(\Omega-\ox)$ is a linear subspace of $\R^n$.\\
{\rm(d)} $\ox\in\Omega$ and $\overline{\cone}(\Omega-\ox)$ is a linear subspace of $\R^n$.\\
{\rm(e)} $\ox\in\Omega$ and the normal cone $N(\ox;\Omega)$ is a subspace of $\R^n$.
\end{Theorem}
{\bf Proof.} (a)$\Longrightarrow$ (b): Suppose that $\ox\in\ri(\Omega)$ and fix $x\in\Omega$ with $x\ne\ox$. It follows from \eqref{ri} that $\ox\in\Omega$ and there exists $\delta>0$ such that
\begin{equation}\label{affc}
\B(\ox;\delta)\cap\aff(\Omega)\subset\Omega.
\end{equation}
Let us choose $0<t<1$ so small that $u:=\ox+t(\ox-x)\in\B(\ox;\delta)$. Hence $u\in\aff(\Omega)$, and we get from \eqref{affc} that $u\in\Omega$. This tells us that
\begin{equation*}
\ox=\frac{t}{1+t}x+\frac{1}{1+t}u\in(x,u),
\end{equation*}
which therefore verifies the conclusions in (b).

(b)$\Longrightarrow$ (c): It suffices to show that for every $a\in K:=\cone(\Omega-\ox)$ we get $-a\in K$. Fix any $a\in K$ and find by definition such $t\ge 0$ and $w\in\Omega$ that $a=t(w-\ox)$. If  $w=\ox$, we have $a=0$ and hence $-a=0\in K$. In the case where $w\ne\ox$, take $u\in\Omega$ with $\ox\in(w,u)$ and then find $\gamma>0$ for which $w=\ox+\gamma(\ox-u)$. It follows therefore that
\begin{equation*}
-a=t(\ox-w)=-t\gamma(\ox-u)=t\gamma(u-\ox)\in K,
\end{equation*}
which implies that $K$ is a linear subspace of $\R^n$.

(c)$\Longrightarrow$ (d): This implication is obvious since every linear subspace is closed in $\R^n$.

(d)$\Longrightarrow$ (e): Fix any $v\in N(\ox;\Omega)$ and show that $-v\in N(\ox;\Omega)$. We have from \eqref{nc} that
\begin{equation*}
\la v,x-\ox\ra\le 0\;\mbox{\rm for all }\;x\in\Omega.
\end{equation*}
Using this and denoting $K:=\overline{\cone}(\Omega-\ox)$ tells us that $\la v,z\ra\le 0$ for all $z\in K$. Since $K$ is a subspace, for any $x\in \Omega$ we get that $\ox-x\in K$, and hence
\begin{equation*}
\ox-x=\lim_{k\to\infty}t_k(w_k-\ox),
\end{equation*}
where $t_k\ge 0$ and $w_k\in\Omega$ for every $k\in\N$. It follows therefore that
\begin{equation*}
\la-v,x-\ox\ra =\lim_{k\to\infty}t_k\la v,w_k-\ox\ra\le 0,
\end{equation*}
which yields $-v\in N(\ox;\Omega)$. Thus $N(\ox;\Omega)$ is a linear subspace in $\R^n$.

(e)$\Longrightarrow$ (a): Assuming that (e) is satisfied and arguing by contradiction, suppose that $\ox\notin\ri(\Omega)$. It follows from Theorem~\ref{pst1} that there exists $v\in\R^n$ such that
\begin{equation}\label{cd1}
\la v,x\ra\le\la v,\ox\ra\;\mbox{\rm for all }\;x\in\Omega,
\end{equation}
and furthermore there exists $\tilde{x}\in\Omega$ for which $\la v,\tilde{x}\ra<\la v,\ox\ra$. Then \eqref{cd1} implies that $v\in N(\ox;\Omega)$, and so $-v\in N(\ox;\Omega)$. This yields
\begin{equation*}
\la-v,\tilde{x}-\ox\ra=\la v,\ox\ra-\la v,\tilde{x}\ra\le 0,
\end{equation*}
which gives us a contradiction, which shows that (a) holds. $\h$\vspace*{0.05in}

The obtained finite-dimensional characterizations of relative interior motivate the major extensions of this notion to infinite dimensions, which are considered in the next section.

\section{Intrinsic Relative and Quasi-Relative Interiors}
\setcounter{equation}{0}

We start with the following basic definitions used throughout the whole paper.

\begin{Definition}{\bf(extended relative interiors in infinite dimensions).}\label{qri} Let $X$ be an LCTV space, and let $\Omega$ be a nonempty convex subset of $X$. Then:\\[1ex]
{\rm(a)} The {\sc intrinsic relative interior} of $\Omega$ is the set
\begin{equation*}
\mbox{\rm iri}(\Omega):=\big\{x\in\Omega\;\big|\;\mbox{\rm cone}(\Omega-x)\;\mbox{\rm is a subspace of }\;X\big\}.
\end{equation*}
{\rm(b)} The {\sc quasi-relative interior} of $\Omega$ is the set
\begin{equation*}
\mbox{\rm qri}(\Omega):=\big\{x\in\Omega\;\big|\;\Bar{\mbox{\rm cone}}(\Omega-x)\;\mbox{\rm is a subspace of }\;X\big\}.
\end{equation*}
{\rm(c)} We say that a convex set $\Omega\subset X$ is {\sc quasi-regular} if $\qri(\Omega)={\rm iri}(\Omega)$.

By convention we set $\mbox{\rm iri}(\emptyset)=\mbox{\rm qri}(\emptyset):=\emptyset$.
\end{Definition}

Due to Theorem~\ref{pts2}, both notions in Definition~\ref{qri}(a,b) reduce to the relative interior of $\Omega$ in finite-dimensional spaces. The one in (a) has been known under the name ``intrinsic core" \cite{hol} (which may be confusing; see \cite{BG}) and also under the name ``pseudo-relative interior" \cite{BG}, which seems to be confusing as well since ``pseudo" means ``false". Following \cite{bao-mor} and \cite{m-book}, we prefer to use the term {\em intrinsic relative
interior} of $\Omega$ and denote this set by $\iri(\Omega)$.

As mentioned in Section~1, the notion of {\em quasi-relative interior} for convex sets was introduced in \cite{bl} and then was studied therein and other publications. It follows directly from the definitions that for any convex set $\Omega\subset X$ we have
\begin{equation*}
\ri(\Omega)\subset\iri(\Omega)\subset\qri(\Omega).
\end{equation*}
In Definition~\ref{qri}(c) we designate the property $\qri(\Omega)={\rm iri}(\Omega)$ by labeling the sets satisfying this condition as {\em quasi-regular} ones. The latter property plays an important role in the subsequent results of the paper. Some sufficient conditions for the quasi-regularity property of convex sets are presented below.

The following example taken from \cite{BG} demonstrates that the sets $\qri(\Omega)$ and $\iri(\Omega)$  may be different in rather simple situations in the framework of Hilbert spaces.

\begin{Example}{\bf(difference between quasi-relative and intrinsic relative interiors).}\label{ri-diff} {\rm Consider the classical space of sequences $X:=\ell^2$ and its convex subset
\begin{equation*}
\Omega:=\Big\{x=(x_k)\in X\;\Big|\;\|x\|_1:=\sum_{k=1}^\infty|x_k|\le 1\Big\}.
\end{equation*}
Then we have $\mbox{\rm iri}(\Omega)=\{x\in X\;|\;\|x\|_1<1\}$, while
\begin{equation*}
\qri(\Omega)=\Omega\setminus\big\{x=(x_k)\in X\;\big|\;\|x\|_1=1,\;x_k=0\;\mbox{ for all }\;k\ge k_0\;\mbox{\rm with some }\;k_0\in\N\big\}.
\end{equation*}}
\end{Example}

The next theorem provides some characterizations of intrinsic relative and quasi-relative interiors of convex sets in LCTV spaces and shows that the equivalences (b)$\Longleftrightarrow$(c) and (d)$\Longleftrightarrow$(e) of Theorem~\ref{ri} hold true in infinite dimensions.

\begin{Theorem}{\bf(characterizations of intrinsic relative and quasi-relative interiors in LCTV space).}\label{Thrpriqri} Let $X$ be an LCTV space, and let $\Omega$ be a nonempty convex subset of $X$ with $\ox\in\Omega$. Then we have the equivalences:\\[1ex]
{\rm(a)} $\ox\in\mbox{\rm iri}(\Omega)$ if and only if for every point $x\in\Omega$ with $x\neq \ox$ there exists another point $u\in\Omega$ such that
$\ox\in(x,u)$.\\
{\rm(b)} $\ox\in\qri(\Omega)$ if and only if the normal cone $N(\ox;\Omega)$ is a subspace of $X^*$.
\end{Theorem}
{\bf Proof.} The equivalence in (a) follows from \cite[Lemma~2.3]{BG}, while the one in (b) is proved in \cite[Proposition~2.8]{bl}. $\h$\vspace*{0.05in}

Our next goal is to provide some sufficient conditions that ensure the quasi-regularity property of convex sets. To proceed, we recall a fundamental notion of general variational analysis in infinite-dimensional spaces that plays a crucial role in furnishing limiting procedures; see \cite{m-book1} for a comprehensive study and applications in convex and nonconvex settings.

\begin{Definition}{\bf(sequential normal compactness of sets).}\label{snc} A nonempty closed convex subset $\Omega$ of a normed space $X$ is called {\sc sequentially normally compact} ${\rm(SNC)}$ {\sc at} $\ox\in\Omega$ if for any sequence $\{x_k\}$ in $\Omega$ converging to $\ox$ we have the implication
\begin{equation*}
\big[x^*_k\in N(x_k;\Omega),\;x^*_k\xrightarrow{w^*}0\big]\Longrightarrow\big[\|x^*_k\|\to 0\big]\;\mbox{ as }\;k\to\infty,
\end{equation*}
where symbol $w^*$ indicates the weak$^*$ topology on $X^*$. The set $\Omega$ is said to be {\sc SNC} if it has this property at any point $\ox\in\Omega$.
\end{Definition}

A (convex) set $\Omega\subset X$ is surely SNC if $\mbox{\rm int}(\Omega)\ne\emp$; see \cite[Proposition~1.25]{m-book1}. If $X$ is Banach and $\mbox{\rm ri}(\Omega)\ne\emp$, we have the following more delicate characterization: $\Omega$ is SNC if and only if the closure of its affine hull $\Bar{\aff}(\Omega)$ is {\em finite-codimensional}; see \cite[Theorem~1.21]{m-book1}. Note that, in the framework of arbitrary normed spaces, a similar characterization is available for a generally more restrictive notion of the {\em compactly epi-Lipschitzian} (CEL) property for closed convex sets $\Omega$ with $\mbox{\rm ri}(\Omega)\ne\emp$: $\Omega$ is CEL if and only if its affine hull is a closed finite-codimensional subspace of $X$; see \cite{blm} and more discussions in \cite[Remark~1.27]{m-book1}.

Prior to deriving sufficient conditions for quasi-regularity, we present three lemmas of their own interest that are used in what follows. The first one is an intrinsic relative interior version in LCTV spaces of the well-known property for relative interiors in finite dimensions.

\begin{Lemma}{\bf(convex sets with nonempty intrinsic relative interiors in LCTV spaces).}\label{lm2separate} Let $X$ be an LCTV space, and let $\Omega\subset X$ be a nonempty closed convex set with $0\in\Omega\setminus\mbox{\rm iri}(\Omega)$. If ${\rm iri}(\Omega)\ne\emp$, then $\overline{\aff}(\Omega)$ is a closed subspace of $X$ and there is a sequence $\{x_k\}\subset -\Omega$ with $x_k\notin\Omega$ and $x_k\to 0$ as $k\to\infty$.
\end{Lemma}
{\bf Proof.} Using ${\rm iri}(\Omega)\ne\emp$ and $0\in\Omega\setminus\mbox{\rm iri}(\Omega)$, let us show that there exists a nonzero vector $x_0\in{\rm iri}(\Omega)$ such that $-tx_0\notin\Omega$ for all $t>0$. Arguing by contradiction, suppose that $-tx_0\in \Omega$ for some $t>0$. Then it follows from \cite[Lemma~3.1]{BG} that
\begin{equation*}
0=\frac{t}{1+t}x_0+\frac{1}{1+t}\big(-tx_0\big)\in\mbox{\rm iri}(\Omega),
\end{equation*}
which clearly contradicts the imposed assumption on $0\notin{\rm iri}(\Omega)$. Denoting now $x_k:=-(x_0/k)\in -\Omega$ tells us that $x_k\notin{\Omega}$ for every $k\in\N$ and that $x_k\to 0$ as $k\to\infty$.$\h$

The second lemma, which is broadly employed in the subsequent material, provides an equivalent description of quasi-relative interiors via {\em proper separation} of convex sets in the general LCTV space framework; see \cite[Theorem~2.3]{Flores1}.

\begin{Lemma}{\bf(proper separation description of quasi-relative interiors in LCTV spaces).}\label{proper separation}
Let $\Omega$ be a convex set in an LCTV space $X$ with $\ox\in\Omega$. Then $\ox\notin\mbox{\rm qri}(\Omega)$ if and only if the sets $\{\ox\}$ and $\Omega$ can be properly separated.
\end{Lemma}
{\bf Proof.} As stated in Theorem~\ref{Thrpriqri}, $\ox\in\mbox{\rm qri}(\Omega)$ if and only if the normal cone $N(\ox;\Omega)$ is a subspace of $X^*$. Thus we get that $\ox\notin\mbox{\rm qri}(\Omega)$ if and only if there exists $x^*\in N(\ox;\Omega)$ with $-x^*\notin N(\ox;\Omega)$. It follows from the normal cone definition \eqref{nc} that $\la x^*,x\ra\le\la x^*,\ox\ra$ for all $x\in\Omega$. Then the inclusion $-x^*\notin N(\ox;\Omega)$ gives us such $x_0\in\Omega$ that $\la-x^*,x_0\ra>\la-x^*,\ox\ra$, which reads as $\la x^*,x_0\ra<\la x^*,\ox\ra$ and hence justifies the claim. $\h$\vspace*{0.05in}

The third lemma provides a useful version of the {\em strict separation} theorem relative to closed subspaces of Hilbert spaces.

\begin{Lemma}{\bf(strict separation relative to subspaces in Hilbert spaces).}\label{Proseparate1} Let $L$ be a closed subspace of a Hilbert space $X$, and let $\Omega\subset L$ be a nonempty convex set with $\bar{x}\in L$ and $\bar{x}\not\in\Bar{\Omega}$. Then there exists a nonzero vector $u\in L$ such that
\begin{equation*}
\sup\big\{\langle u,x\rangle\;\big|\;x\in\Omega\big\}<\langle u,\bar{x}\rangle.
\end{equation*}
\end{Lemma}
{\bf Proof.} Since $\ox\notin\Bar{\Omega}$, we see that the sets $\{\ox\}$ and $\Omega$ are strictly separated in $X$, which means that there exists a vector $h\in X$ such that
\begin{equation}\label{eqlm1}
\sup\big\{\la h,x\ra\;\big|\;x\in\Omega\big\}<\la h,\ox\ra.
\end{equation}
It is well known that the space $X$ can be represented as the direct sum $X=L\oplus L^{\bot}$, where
\begin{equation*}
L^{\bot}:=\big\{w\in X\;\big|\;\la w,x\ra=0\;\text{ for all }\;x\in L\big\}.
\end{equation*}
If $h\in L^{\bot}$, then \eqref{eqlm1} immediately gives us a contradiction. Thus $h\in X$ can be represented as $h=u+w$ with $0\ne u\in L$ and $w\in L^{\bot}$. This implies that for each $x\in\Omega\subset L$ we have
\begin{equation*}
\begin{aligned}
\la u,x\ra &=\la u,x\ra+\la w,x\ra=\la h,x\ra\le\sup\big\{\la h,x\ra\;\big|\;x\in\Omega\big\}<\la h,\ox\ra\\
&=\la u+w,\ox\ra=\la u,\ox\ra,
\end{aligned}
\end{equation*}
which shows that $\sup\{\la u,x\ra\;|\;x\in\Omega\}<\la u,\ox\ra$ with $u\ne 0$. $\h$\vspace*{0.05in}

Now we are ready to present sufficient conditions for quasi-regularity of convex sets.\vspace*{-0.05in}

\begin{Theorem}{\bf(sufficient conditions for quasi-regularity).}\label{qreg} Let $\Omega\subset X$ be a nonempty convex set. Then it is quasi-regular under each of the following conditions:\\[1ex]
{\rm (a)} $X$ is a finite-dimensional space.\\
{\rm (b)} $X$ is an LCTV space and $\mbox{\rm int}(\Omega)\ne\emp$.\\
{\rm (c)} $X$ is an LCTV space and $\mbox{\rm ri}(\Omega)\ne\emp$.\\
{\rm (d)} $X$ is a Hilbert space, $\Omega$ is closed and SNC with $\iri(\Omega)\ne\emp$.
\end{Theorem}
{\bf Proof.} The quasi-regularity of $\Omega$ in case (a) is proved in Theorem~\ref{pts2}. The fulfillment of this property under (b) and (c) follows from \cite[Theorem~2.12]{BG}. Let us now show that $\Omega$ is quasi-regular under the assumptions in (d). We split the proof into the following two steps.

{\bf Step~1.} First we verify that in the case where $0\notin{\rm iri}(\Omega)$, the sets $\Omega$ and $\{0\}$ are properly separated, i.e., there is a nonzero vector $a\in X$ such that
\begin{equation}\label{pro-sep}
\sup\big\{\langle a,x\rangle\;\big|\;x\in\Omega\big\}\le 0\;\mbox{ and }\;\inf\big\{\langle a,x\rangle\;\big|\;x\in\Omega\big\}<0,
\end{equation}
which is equivalent therefore to $0\notin\qri(\Omega)$ (assuming $0\in \Omega$) by Lemma~\ref{proper separation}.

If $0\not\in\Omega$, this statement is trivial. Suppose now that $0\in\Omega\setminus\text{\rm iri}(\Omega)$. Letting $L:=\overline{\mbox{\rm aff}}(\Omega)$ and employing Lemma~\ref{lm2separate} tell us that $L$ is a subspace of $X$ and that there exists a sequence $\{x_k\}\subset L$ for which $x_k\notin {\Omega}$ and $x_k\to 0$ as $k\to\infty$. By Lemma~\ref{Proseparate1} we find a sequence $\{v_k\}\subset L$ with $v_k\ne 0$ such that
\begin{equation*}
\sup\big\{\langle v_k,x\rangle\;\big|\;x\in\Omega\big\}<\la v_k,x_k\ra\;\mbox{ whenever }\;k\in\N.
\end{equation*}
Denote $w_k:=\frac{v_k}{\Vert v_k\Vert}\in L$ with $\|w_k\|=1$ as $k\in\N$ and observe that
\begin{equation}\label{eq1Lm1.5wk}
\langle w_k,x\rangle<\la w_k,x_k\ra\le\|w_k\|\cdot\|x_k\|=\varepsilon_k\;\mbox{\rm for all }\;x\in\Omega,
\end{equation}
where $\varepsilon_k:=\|x_k\|\downarrow 0$ as $k\to\infty$. Since $\{w_k\}$ is bounded, we let $k\to\infty$ in (\ref{eq1Lm1.5wk}) and suppose without loss of generality that $w_k\xrightarrow{w}a\in L$, where the symbol ``$w$" indicates the weak convergence in the Hilbert space $X$. This implies that
\begin{equation}\label{eq2Lm1.5}
\sup\big\{\langle a,x\rangle\;\big|\;x\in\Omega\big\}\le 0.
\end{equation}

To verify further the strict inequality
\begin{equation*}
\inf\big\{\langle a,x\rangle\;\big|\;x\in\Omega\big\}<0,
\end{equation*}
it suffices to show that there exists $\ox\in\Omega$ with $\langle a,\ox\rangle<0$. Arguing by contradiction, suppose that $\langle a,x\rangle\ge 0$ for all $x\in\Omega$ and deduce from \eqref{eq2Lm1.5} that $\langle a,x\rangle=0$ whenever $x\in\Omega$. Since $a\in L=\overline{\aff}(\Omega)$, there exists a sequence $a_j\to a$ as $j\to\infty$ with $a_j\in\aff(\Omega)$. The latter inclusion can be equivalently rewritten as
\begin{equation*}
a_j=\sum_{i=1}^{m_j}\lambda^j_i\omega^j_i\;\mbox{ with }\;\sum_{i=1}^{m_j}\lambda^j_i=1\;\mbox{ and }\;\omega^j_i\in\Omega\;\mbox{ for }\;i=1,\ldots,m_j,
\end{equation*}
which readily implies the equalities
\begin{equation*}
\langle a, a_j\rangle=\sum_{i=1}^{m_j}\lambda^j_i\langle a,\omega^j_i\rangle=0.
\end{equation*}
The passage to the limit as $j\to\infty$ gives us $\|a\|^2=0$, and so $a=0$.

Now we apply the result of \cite[Theorem~3.1.2]{z}, which ensures by using \eqref{eq1Lm1.5wk} the existence of $b_k\in\Omega$ and $u_k\in X$ such that
\begin{equation}\label{eq3Lm1.5}
u_k\in N(b_k;\Omega),\;\|b_k\|\le\sqrt{\varepsilon_k},\;\mbox{ and }\;\|u_k-w_k\|\le\sqrt{\varepsilon_k}.
\end{equation}
Since $\|w_k\|=1$, it follows from \eqref{eq3Lm1.5} that $\|u_k\|\to 1$. Furthermore, we get from $w_k\xrightarrow{w}0$, $\varepsilon_k\downarrow 0$, and \eqref{eq3Lm1.5} that $u_k\xrightarrow{w}0$ as $k\to\infty$. Remembering that $\Omega$ enjoys the SNC property, it follows from \eqref{eq3Lm1.5} that $\|u_k\|\to 0$, which clearly contradicts the condition $\|u_k\|\to 1$ as $k\to\infty$. This tells us that there exists $\ox\in\Omega$ such that $\la a,\ox\ra<0$, and hence $a\ne 0$. It justifies the proper separation of $\Omega$ and $\{0\}$ in \eqref{pro-sep} and shows therefore that $0\notin\qri(\Omega)$.\\[1ex]
{\bf Step 2.} To verify the quasi-regularity of $\Omega$, we need to prove that $\mbox{\rm qri}(\Omega)\subset\mbox{\rm iri}(\Omega)$. Picking any $\ox\in \mbox{\rm qri}(\Omega)$ gives us $0\in\qri(\Omega-\ox)$. Since $\Omega$ is SNC, this property holds for $\Omega-\ox$ as well. It allows us to deduce from the proof in Step~1 that $0\in\mbox{\rm iri}(\Omega-\ox)$, and so $\ox\in\mbox{\rm iri}(\Omega)$. This shows that $\mbox{\rm qri}(\Omega)\subset\mbox{\rm iri}(\Omega)$ and thus completes the proof of the theorem. $\h$

\section{Quasi-Relative Interiors of Convex Graphs}
\setcounter{equation}{0}

In this section we start studying properties of set-valued mappings with convex graphs that involve the (extended) relative interior constructions defined above. Remember that a set-valued mapping $F\colon X\tto Y$ between LCTV spaces is associated with its {\em graph}
\begin{equation*}
\gph(F):=\big\{(x,y)\in X\times Y\;\big|\;y\in F(x)\big\},
\end{equation*}
and it is called {\em convex} if its graph is a convex subset of the product space $X\times Y$. We also consider the {\em domain} of $F$ defined by
\begin{equation*}
\dom(F):=\big\{x\in X\;\big|\;F(x)\ne\emp\big\}.
\end{equation*}

Let us first present a new proof, based on proper separation, of the important result by Rockafellar \cite[Theorem~6.8]{r} in finite dimensions (see also \cite[Proposition~2.43]{rw} for yet another proof), which can be formulated in the following way.

\begin{Theorem}{\bf(Rockafellar's theorem on relative interiors of convex graphs).}\label{TheoRoc1} Let $F\colon\R^m\tto\R^n$ be a convex set-valued mapping. Then we have the representation
\begin{equation}\label{roc}
\ri\big(\gph(F)\big)=\big\{(x,y)\in\R^m\times\R^n\;\big|\;x\in\ri\big(\dom(F)\big),\;y\in\ri\big(F(x)\big)\big\}.
\end{equation}
\end{Theorem}
{\bf Proof.} We first prove the inclusion ``$\subset$" in \eqref{roc}. Picking $(\ox,\oy)\in\ri(\gph(F))$ and arguing by contradiction, suppose that $\ox\notin\ri(\dom (F))$. It follows from the relative interior descriptions of Theorem~\ref{pst1} that there exists $v\in\mathbb\R^m$ such that
\begin{equation*}
\la v,x\ra\le\la v,\ox\ra\;\mbox{ for all }\;x\in\dom(F)
\end{equation*}
together with a vector $x_0\in\dom(F)$ satisfying
\begin{equation*}
\la v,x_0\ra<\la v,\ox\ra.
\end{equation*}
Then for all $(x,y)\in\gph(F)$ we have $x\in\dom(F)$ and
\begin{equation*}
\la(v,0),(x,y)\ra=\la v,x\ra\le\la v,\ox\ra=\la(v,0),(\ox,\oy)\ra.
\end{equation*}
Taking now any fixed $y_0\in F(x_0)$ gives us
\begin{equation*}
\la(v,0),(x_0,y_0)\ra=\la v,x_0\ra<\la v,\ox\ra=\la(v,0),(\ox,\oy)\ra.
\end{equation*}
This shows that the sets $\gph(F)$ and $\{(\ox,\oy)\}$ are properly separated, and hence $(\ox,\oy)\notin\ri(\gph(F))$ by Theorem~\ref{pst1}. The obtained contradiction verifies that $\ox\in\ri(\dom(F))$.

To show now that $\oy\in\ri(F(\ox))$, we deduce from $(\ox,\oy)\in\ri(\gph(F))$  due to definition \eqref{ri} that there exists $\delta>0$ such that
\begin{equation*}
\big[\B(\ox;\delta)\times\B(\oy;\delta)\big]\cap\aff\big(\gph(F)\big)\subset\gph(F),
\end{equation*}
which yields in turn the inclusion
\begin{equation}\label{i1}
\big[\{\ox\}\times\B(\oy;\delta)\big]\cap\aff\big(\gph(F)\big)\subset\gph(F).
\end{equation}
It follows furthermore from the definition of the affine hull that
\begin{equation*}
\{\ox\}\times\aff\big(F(\ox)\big)\subset\aff\big(\gph(F)\big),
\end{equation*}
which being combined with \eqref{i1} leads us to the relationships
\begin{align*}
\{\ox\}\times\big[(\B(\oy;\delta)\cap\aff\big(F(\ox)\big)\big]&=\big[\{\ox\}\times\B(\oy;\delta)\big]\cap\big[\{\ox\}\times\aff\big(F(\ox)\big)\big]\\
&\subset\big[\{\ox\}\times\B(\oy;\delta)\big]\times\aff\big(\gph(F)\big)\subset\gph(F).
\end{align*}
This shows that $\B(\oy,\delta)\cap\aff(F(\ox))\subset F(\ox)$ and so verifies that $\oy\in\ri(F(\ox))$.

To prove next the opposite inclusion in \eqref{roc}, fix $\ox\in\ri(\dom(F))$ and $\oy\in\ri(F(\ox))$. Arguing by contradiction, suppose that $(\ox,\oy)\notin\ri(\gph(F))$ and then find by Theorem~\ref{pst1} such a pair $(u,v)\in\R^m\times\R^n$ that
\begin{equation}\label{1}
\la u,x\ra+\la v,y\ra\le\la u,\ox\ra+\la v,\oy\ra\;\mbox{ whenever }\;y\in F(x).
\end{equation}
In addition, it follows from the proper separation that there is $(x_0,y_0)\in\gph(F)$ satisfying
\begin{equation}\label{2}
\la u,x_0\ra+\la v,y_0\ra<\la u,\ox\ra+\la v,\oy\ra.
\end{equation}
Letting $x=\ox$ in \eqref{1} yields $\la v,y\ra\le\la v,\oy\ra$ for all $y\in F(\ox)$.  Since $\ox\in\ri(\dom(F))$ and $x_0\in \dom(F)$, we deduce from Theorem~\ref{pts2}(a) that there exists $\tilde{x}\in\dom(F)$ such that $\ox=tx_0+(1-t)\tilde{x}$ for some $t\in(0,1)$. Choose $\tilde{y}\in F(\tilde{x})$ and consider the convex combination
\begin{equation*}
y^\prime:=ty_0+(1-t)\tilde{y},
\end{equation*}
where $y^\prime\in F(\ox)$ since $\gph(F)$ is convex. Since $\tilde{x}\in\dom(F)$ we use \eqref{1} and \eqref{2} to get
\begin{align*}
&\la u,\tilde{x}\ra+\la v,\tilde{y}\ra\le\la u,\ox\ra+\la v,\oy\ra,\\
&\la u,x_0\ra+\la v,y_0\ra<\la u,\ox\ra+\la v,\oy\ra.
\end{align*}
Multiplying the first inequality above by $1-t$ and the second one by $t$, and then adding them together gives us the condition
\begin{equation*}
\la u,\ox\ra+\la v,y^\prime\ra<\la u,\ox\ra+\la v,\oy\ra,
\end{equation*}
which yields $\la v,y^\prime\ra<\la v,\oy\ra$. Thus the sets $\{\oy\}$ and $F(\ox)$ are properly separated. Applying Theorem~\ref{pst1} tells us that $\oy\notin\ri(F(\ox))$, a contradiction that verifies $(\ox,\oy)\in\ri(\gph(F))$. $\h$

As a direct consequence of Theorem~\ref{TheoRoc1}, we get the following useful result on calculating relative interiors for epigraphs of extended-real-valued convex functions on $\R^n$. Recall that the {\em domain} and {\em epigraph} of a function $f\colon X\to\oR:=(-\infty,\infty]$ are defined, respectively, by
\begin{equation*}
\dom(f):=\big\{x\in X\;\big|\;f(x)<\infty\big\}\;\mbox{ and }\;\epi(f):=\big\{(x,\al)\in X\times\R\;\big|\;\al\ge f(x.)\big\}
\end{equation*}
The convexity of $f$ corresponds to the convexity of $\epi(f)$, and $f$ is {\em proper} if $\dom(f)\ne\emp$.

\begin{Corollary}{\bf(relative interiors of epigraphs for convex functions on $\R^n$).}\label{epi-roc} Let $f\colon\R^n\to\oR$  be a convex function. Then we have
\begin{equation}\label{epi}
\ri\big(\epi(f)\big)=\big\{(x,\lambda)\in\R^{n+1}\;\big|\;x\in\ri(\dom(f)),\;\lambda>f(x)\big\}.
\end{equation}
\end{Corollary}
{\bf Proof.} Given $f\colon\R^n\to\oR$, define the set-valued mapping $F\colon\R^n\tto\R$ by $F(x):=\big[f(x),\infty\big)$ for all $x\in\R^n$. Then $\dom(F)=\dom(f)$, $\gph(F)=\epi(f)$, and for any $x\in\dom(f)$ we have $\ri(F(x))=\big(f(x),\infty\big)$. Then \eqref{epi} follows directly from Theorem~\ref{TheoRoc1}. $\h$\vspace*{0.05in}

Next we extend Rockafellar's theorem to the general setting of LCTV spaces. As the reader can see, the extension below, which formulated in terms of the quasi-relative interior, is more involved in infinite dimensions while giving us two inclusion counterparts of \eqref{roc} under different quasi-regularity assumptions. Our proof is based on the proper separation description of the quasi-relative interior established in Lemma~\ref{proper separation}.

\begin{Theorem}{\bf(quasi-relative interiors of convex graphs in LCTV space).}\label{rocthm2} Let $F\colon X\tto Y$ be a convex set-valued mapping between LCTV spaces. The following hold:\\[1ex]
{\rm(a)} If the graph $\gph(F)$ is quasi-regular, then we have the inclusion
\begin{equation*}
\mbox{\rm qri}\big(\gph(F)\big)\subset\big\{(x,y)\in X\times Y\;\big|\;x\in\mbox{\rm qri}\big(\dom(F)\big),\;y\in\mbox{\rm qri}\big(F(x)\big)\big\}.
\end{equation*}
{\rm(b)} If the domain $\dom(F)$ is quasi-regular, then we have the opposite inclusion
\begin{equation*}
\mbox{\rm qri}\big(\gph(F)\big)\supset\big\{(x, y)\in X\times Y\;\big|\;x\in\mbox{\rm qri}\big(\dom(F)\big),\;y\in\mbox{\rm qri}\big(F(x)\big)\big\}.
\end{equation*}
\end{Theorem}
{\bf Proof.} To verity the inclusion ``$\subset$" in (a), pick $(\ox,\oy)\in\qri(\gph(F))$ and suppose on the contrary that $\ox\notin\qri(\dom(F))$. By Lemma~\ref{proper separation} we find $v^*\in X^*$ such that
\begin{equation*}
\la v^*,x\ra\le\la v^*,\ox\ra\;\mbox{ whenever }\;x\in\dom(F)
\end{equation*}
and also have $x_0\in\dom(F)$ for which the strict inequality
\begin{equation*}
\la v^*,x_0\ra<\la v^*,\ox\ra
\end{equation*}
is satisfied. Then for all $(x,y)\in\gph(F)$ we get
\begin{equation*}
\la(v^*,0),(x,y)\ra=\la v^*,x\ra\le\la v^*,\ox\ra=\la(v^*,0),(\ox,\oy)\ra,
\end{equation*}
and for each $y_0\in F(x_0)$ arrive at the conditions
\begin{equation*}
\la(v^*,0),(x_0,y_0)\ra=\la v^*,x_0\ra<\la v^*,\ox\ra=\la (v^*,0),(\ox,\oy)\ra.
\end{equation*}
This shows that the sets $\gph(F)$ and $\{(\ox,\oy)\}$ are properly separated, and hence $(\ox,\oy)\notin\qri\big(\gph(F)\big)$ by Lemma~\ref{proper separation}. The obtained contradiction tells us that $\ox\in \qri(\dom(F))$.

To proceed further with the proof of (a), let us verify that $\oy\in\mbox{\rm qri}(F(\ox))$. Fix any $y\in F(\ox)$ with $y\ne\oy$, and so $(\ox,y)\in\gph(F)$. The assumed quasi-regularity of $\gph(F)$ gives us $(\tilde{x},\tilde{y})\in\gph(F)$ and $t\in(0,1)$ such that
\begin{equation*}
(\ox,\oy)=t(\ox,y)+(1-t)(\tilde{x},\tilde{y}).
\end{equation*}
This yields $\tilde{x}=\ox$ and $\oy=ty+(1-t)\tilde{y}$ with $\tilde{y}\in F(\ox)$. It follows therefore that $\oy\in\mbox{\rm iri}(F(\ox))\subset\mbox{\rm qri}(F(\ox))$, which completes the proof of (a).

To verify now assertion (b) under the quasi-regularity of $\dom(F)$, fix $\ox\in\qri(\dom(F))$ and $\oy\in\qri\big(F(\ox)\big)$. Arguing by contradiction, suppose that $(\ox,\oy)\notin\qri(\gph(F))$. Then Lemma~\ref{proper separation} ensures the existence of $(u^*,v^*)\in X^*\times Y^*$ such that
\begin{equation}\label{equaRck2}
\la u^*,x\ra+\la v^*,y\ra\le\la u^*,\ox\ra+\la v^*,\oy\ra\;\mbox{ whenever }\;y\in F(x)
\end{equation}
and also the existence of $(x_0,y_0)\in\gph(F)$ for which
\begin{equation*}
\la u^*,x_0\ra+\la v^*,y_0\ra<\la u^*,\ox\ra+\la v^*,\oy\ra.
\end{equation*}
Letting $x=\ox$ in \eqref{equaRck2} yields $\la v^*,y\ra\le\la v^*,\oy\ra$ for all $y\in F(\ox)$. Using $\ox\in\qri\big(\dom(F)\big)$, $x_0\in\dom(F)$, and the assumed quasi-regularity of $\dom(F)$ allows us to deduce from Theorem~\ref{Thrpriqri}(a) the existence of $\tilde{x}\in\dom(F)$ that ensures the representation $\ox=tx_0+(1-t)\tilde{x}$ with some $t\in(0, 1)$. Pick $\tilde{y}\in F(\tilde{x})$ and define
\begin{equation*}
y^\prime:=ty_0+(1-t)\tilde{y}.
\end{equation*}
Then $y^\prime\in F(\ox)$ by the convexity of $\gph(F)$, and we get
\begin{align*}
&\la u^*,\tilde{x}\ra+\la v^*,\tilde{y}\ra\le\la u^*,\ox\ra+\la v^*,\oy\ra,\\
&\la u^*,x_0\ra+\la v^*,y_0\ra<\la u^*,\ox\ra+\la v^*,\oy\ra.
\end{align*}
Multiply the first inequality above by $1-t$, the second inequality by $t$, and add them together to arrive at the condition
\begin{equation*}
\la u^*,\ox\ra+\la v^*,y^\prime\ra<\la u^*,\ox\ra+\la v^*,\oy\ra,
\end{equation*}
which gives us $\la v^*,y^\prime\ra<\la v^*,\oy\ra$. Thus we conclude that the sets $\{\oy\}$ and $F(\ox)$ are properly separated, and so $\oy\notin\qri(F(\ox))$ by Lemma~\ref{proper separation}. This contradiction shows that $(\ox,\oy)\in\qri(\gph(F))$ and hence completes the proof of the theorem. $\h$

Based on Theorem~\ref{Lemalinearmap}, we arrive now at  the precise calculation of the quasi-relative interiors of convex graphs for set-valued mappings $F\colon X\tto Y$ between LCTV spaces.

\begin{Corollary}{\bf(precise calculation of quasi-relative interiors of convex graphs in LCTV spaces).}\label{qri-precise} Let $F\colon X\tto Y$ be a convex set-valued mapping between LCTV spaces. If both $\gph(F)$ and $\dom(F)$ are quasi-regular, then we have
\begin{equation*}
\mbox{\rm qri}\big(\gph(F)\big)=\big\{(x,y)\in X\times Y\;\big|\;x\in\mbox{\rm qri}\big(\dom(F)\big)\;y\in\mbox{\rm qri}\big(F(x)\big)\big\}.
\end{equation*}
\end{Corollary}

\section{Further Properties of Extended Relative Interiors}
\setcounter{equation}{0}

In this section we continue our study of important properties of the quasi-relative and intrinsic relative interiors of convex set in LCTV spaces, with a particular emphasis on graphical, domain, and epigraphical sets for convex set-valued mappings and extended-real-valued functions. Note that the assumption $\qri(\Omega)\ne\emp$, which is imposed in many results presented in this section, automatically holds for closed convex sets in any separable Banach spaces due to the fundamental result of \cite[Theorem~2.19]{bl}.

The first statement describes behavior of quasi-relative interiors of convex sets under applying linear continuous operators.\vspace*{-0.1in}

\begin{Theorem}{\bf(quasi-relative interiors under linear mappings).}\label{Lemalinearmap} Let $A\colon X\to Y$ be a linear continuous mapping between two LCTV spaces, and let $\Omega\subset X$ be a convex such that $\iri(\Omega)\neq\emptyset$ and $A(\Omega)$ is quasi-regular. Then we have
\begin{equation*}
A\big(\qri(\Omega)\big)=\qri\big(A(\Omega)\big).
\end{equation*}
\end{Theorem}
{\bf Proof.} The inclusion $A(\qri(\Omega))\subset\qri(A(\Omega))$ is known; see \cite[Proposition~2.21]{bl}. To verify the opposite inclusion, we use \cite[Proposition~6.3.2(i)]{ktz} and get
\begin{equation*}
\qri(A(\Omega))=\iri(A(\Omega))=A(\iri(\Omega))\subset A(\qri(\Omega)),
\end{equation*}
which completes the proof.  $\h$\vspace*{0.05in}

The next corollary on representing quasi-relative interiors of set differences is a direct consequence of Theorem~\ref{Lemalinearmap}.\vspace*{-0.05in}

\begin{Corollary}{\bf(quasi-relative interiors of set differences).}\label{Corrosub} Let $\Omega_1$ and $\Omega_2$ be convex subsets of an LCTV space $X$. Suppose that $\iri(\Omega_1)\ne\emp$, $\iri(\Omega_2)\ne\emp$, and $\Omega_1-\Omega_2$ is quasi-regular. Then we have the equality
\begin{equation*}
\qri(\Omega_1-\Omega_2)=\qri(\Omega_1)-\qri(\Omega_2).
\end{equation*}
\end{Corollary}
{\bf Proof.} Define the linear continuous mapping $A\colon X\times X\to X$ by $A(x,y):=x-y$, and let $\Omega:=\Omega_1\times\Omega_2$. Then $\mbox{\rm iri}(\Omega)=\mbox{\rm iri}(\Omega_1)\times\mbox{\rm iri}(\Omega_2)\ne\emp$; see \cite[Proposition~2.5]{bl}. Applying Theorem~\ref{Lemalinearmap} gives us
\begin{equation*}
\mbox{\rm qri}(\Omega_1-\Omega_2)=A\big(\qri(\Omega)\big)=\qri\big(A(\Omega)\big)=\qri(\Omega_1)-\mbox{\rm qri}(\Omega_2),
\end{equation*}
and thus we arrive at the claimed equality. $\h$\vspace*{0.05in}

\begin{Theorem}{\bf(proper separation in LCTV spaces via quasi-relative interiors).}\label{Theoqricap} Let $\Omega_1$ and $\Omega_2$ be convex subsets of an LCTV space $X$ such that $\Omega_1\cap\Omega_2\neq\emptyset$. Assume that $\iri(\Omega_1)\ne\emp$, $\iri(\Omega_2)\ne\emp$, and $\Omega_1-\Omega_2$ is quasi-regular. Then $\Omega_1$ and $\Omega_2$ are properly separated if and only if
\begin{equation}\label{qri-sep}
\qri(\Omega_1)\cap\qri(\Omega_2)=\emp.
\end{equation}
\end{Theorem}
{\bf Proof.} Define $\Omega:=\Omega_1-\Omega_2$ and get from Corollary~\ref{Corrosub} that condition \eqref{qri-sep} reduces to
\begin{equation*}
0\notin\qri(\Omega_1-\Omega_2)=\qri(\Omega_1)-\qri(\Omega_2).
\end{equation*}
If \eqref{qri-sep} holds, then $0\notin\qri(\Omega_1-\Omega_2)=\qri(\Omega)$. Since $\Omega_1\cap \Omega_2\neq\emptyset$, we see that $0\in \Omega_1-\Omega_2$, so Lemma~\ref{proper separation} tells us that the sets $\Omega$ and $\{0\}$ are properly separated, which clearly ensures the proper separation of the sets $\Omega_1$ and $\Omega_2$.

To verify the opposite implication, suppose that $\Omega_1$ and $\Omega_2$ are properly separated, which implies that the sets $\Omega=\Omega_1-\Omega_2$ and $\{0\}$ are properly separated as well. Then using Lemma~\ref{proper separation} and Corollary~\ref{Corrosub} tells us that
\begin{equation*}
0\notin\qri(\Omega)=\qri(\Omega_1-\Omega_2)=\qri(\Omega_1)-\qri(\Omega_2),
\end{equation*}
and thus $\qri(\Omega_1)\cap\qri(\Omega_2)=\emp$, which completes the proof. $\h$\vspace*{0.05in}

\begin{Theorem}{\bf(intrinsic relative interiors of convex epigraphs in LCTV spaces).}\label{LMregularepi} Let $f\colon X\to\oR$ be a convex function, where $X$ is an LCTV space. Then we  have the representation
\begin{equation}\label{priepi}
{\rm iri}\big(\epi(f)\big)=\big\{(x,\lambda)\in X\times\R\;\big|\;x\in{\rm iri}\big(\dom(f)\big),\;\lambda>f(x)\big\}.
\end{equation}
\end{Theorem}
{\bf Proof.}  Denoting by $\Omega$ the set on the right-hand side of \eqref{priepi}, we are going to show that
\begin{equation}\label{eq1LMquariepi}
{\rm iri}\big(\epi(f)\big)=\Omega.
\end{equation}
Let us first verify the inclusion ``$\subset$'' in \eqref{eq1LMquariepi}. Pick any $(\ox,\bar{\lambda})\in{\rm iri}(\epi(f))$ and check that $\ox\in{\rm iri}\big(\dom(f)\big)$. Fixing $x\in\dom(f)$ with $x\ne\ox$, we get $(x,\lambda)\in\epi(f)$, where $\lambda:=f(x)$. Then Theorem~\ref{Thrpriqri} ensures the existence of $(u,\gamma)\in\epi(f)$ such that
\begin{equation*}
(\ox,\bar{\lambda})\in\big((x,\lambda),(u,\gamma)\big),
\end{equation*}
which gives us $\ox\in(x,u)$. Applying Theorem~\ref{Thrpriqri} again yields $\ox\in{\rm iri}\big(\dom(f)\big)$.

Let us now show that $f(\ox)<\bar{\lambda}$. Arguing by contradiction, suppose that $\bar{\lambda}=f(\ox)$ and take any $(\ox,\tilde{\lambda})\in\epi(f)$ with $\tilde{\lambda}>f(\ox)$. It means that $(\ox,\tilde{\lambda})\ne(\ox,\bar{\lambda})=\big(\ox,f(\ox)\big)$. Then it follows from Theorem~\ref{Thrpriqri} that there exists $(\bar{u},\bar{\gamma})\in\epi(f)$ with $\big(\ox,f(\ox)\big)\in\big((\ox,\tilde{\lambda}),(\bar{u},\bar{\gamma})\big)$, and hence we can find  $t_0\in(0,1)$ such that
\begin{equation*}
\ox=t_0\ox+(1-t_0)\bar{u}\;\mbox{ and }\;\bar{\lambda}=t_0\tilde{\lambda}+(1-t_0)\bar{\gamma}.
\end{equation*}
Employing the convexity of $f$ shows that
\begin{equation*}
t_0\tilde{\lambda}+(1-t_0)\bar{\gamma}=\bar{\lambda}=f(\ox)\le t_0f(\ox)+(1-t_0)f(\bar{u})<t_0\tilde{\lambda}+(1-t_0)f(\bar{u})
\end{equation*}
thus verifying that $\bar{\gamma}<f(\bar{u})$, and hence $(\bar{u},\bar{\gamma})\notin\epi(f)$. The obtained contradiction tells us that $\overline{\lambda}>f(\ox)$ and therefore justifies the inclusion ${\rm iri}(\epi(f))\subset\Omega$.

To prove next the opposite inclusion in \eqref{eq1LMquariepi}, fix any $(\ox,\bar{\lambda})\in\Omega$ giving us $\ox\in{\rm iri}(\dom(f))$ and $\bar{\lambda}>f(\ox)$. Picking now any $(x,\lambda)\in\epi(f)$ with $(x,\lambda)\ne(\ox,\bar{\lambda})$, let us verify the existence of $(\oy,\bar{\beta})\in\epi(f)$ for which
\begin{equation*}
(\ox,\bar{\lambda})\in\big((x,\lambda),(\oy,\bar{\beta})\big).
\end{equation*}
To proceed, we consider following two cases:\\[1ex]
{\bf Case 1: $x\ne\ox$}. Since $\ox\in{\rm iri}(\dom(f))$ and $\ox\ne x\in\dom(f)$, there exists $u\in\dom(f)$ such that $\ox\in(x,u)$. Choose $\gamma\in\R$ satisfying
\begin{equation*}
(\ox,\bar{\lambda})\in\big((u,\gamma),(x,\lambda)\big)
\end{equation*}
and check that there exists $(\oy,\bar{\beta})\in((u,\gamma),(\ox,\bar{\lambda}))$ with $(\oy,\bar{\beta})\in\epi(f)$. Arguing by contradiction, suppose that for every $(y,\beta)\in((u,\gamma),(\ox,\bar{\lambda}))$ we have $(y,\beta)\notin\epi(f)$. Fix any $t\in(0,1)$ and define the $t$-dependent elements
\begin{equation*}
y_t:=tu+(1-t)\ox\;\mbox{ and }\;\beta_t:=t\gamma+(1-t)\bar{\lambda},
\end{equation*}
for which we get $(y_t,\beta_t)\in((u,\gamma),(\ox,\bar{\lambda}))$. The convexity of $f$ ensures that
\begin{align*}
t\gamma+(1-t)\bar{\lambda}=\beta_t<f(y_t)\le tf(u)+(1-t)f(\ox)\le t f(u)+(1-t)\bar{\lambda}.
\end{align*}
Letting there $t\downarrow 0$ shows that $\bar{\lambda}=f(\ox)$, a contradiction that verifies hence the existence of a pair $(\oy,\bar{\beta})$ with $(\oy,\bar{\beta})\in((u,\gamma),(\ox,\bar{\lambda}))$ and $(\oy,\bar{\beta})\in\epi(f)$. It tells us that $(\ox,\bar{\lambda})\in((x,\lambda),(\oy,\bar{\beta}))$, and so it follows from Theorem~\ref{Thrpriqri} that $(\ox,\bar{\lambda})\in{\rm iri}(\epi(f))$.\\[1ex]
{\bf Case~2: $x=\ox$}. Then we have $\lambda\ne\bar{\lambda}$, and it follows from $\bar{\lambda}>f(\ox)$ that there exists $(\ox,\tilde{\lambda})\in\epi(f)$ with $(\ox,\bar{\lambda})\in((\ox,\lambda),(\ox,\tilde{\lambda}))$. This verifies the representation in \eqref{priepi}. $\h$

\begin{Proposition}\label{affine reg} Let $\Omega$ be a convex set in an LCTV space $X$. Given $x^*\in X^*$ and $b\in \R$, define the extended-real-valued function
\begin{equation*}
 f(x):=\left\{\begin{aligned}
 &\la x^*,x\ra +b\; & \text{if} &\; x\in\Omega,\\
 &\infty \; &  \text{if}\;  &x\notin\Omega.
 \end{aligned}\right.
\end{equation*}
Then we have the representation
\begin{equation}\label{qrinotre}
\qri(\epi(f))= \{(x,\lambda)\in X\times \mathbb{R}\ |\ x\in\qri(\Omega), \lambda>f(x)\}.
\end{equation}

\end{Proposition}
{\bf Proof.} By \eqref{qriepi1} in Theorem \ref{quasi epi}, it suffices to prove the inclusion ``$\subset$'' in \eqref{qrinotre}. Taking $(x_0,\lambda_0)\in \qri(\epi(f)),$ we will show $(x_0,\lambda_0)$ belongs to the set on the right-hand side of \eqref{qrinotre}.
 Define  $F\colon X\tto\mathbb{R}$  by $F(x):=[f(x),\infty)$ for $x\in X$.  Then $\dom(F)=\dom(f)=\Omega,$ and $\gph(F)=\epi(f).$
Following the  proof of part (a) in Theorem~\ref{rocthm2} without assuming the quasi-regularity of $\gph(F)$) yields $x_0\in\qri(\dom(F))=\qri(\Omega)$.  It remains to show that $\lambda_0>f(x_0).$ By contradiction, suppose that $\lambda_0\leq f(x_0)$.  Since $(x_0,\lambda_0)\in\epi(f),$ we have $\lambda_0=f(x_0).$ Then $(x_0,f(x_0))\in\qri(\epi(f))$ and it follows from the definition that $A:=\overline{\cone}\big(\epi(f)-(x_0,f(x_0))\big)$ is a subspace of $X$. Let $a=(x_0,f(x_0)+2)-(x_0,f(x_0))=(0,2)\in A.$ Then $-a=(0,-2)\in A.$ Therefore, there exists a net $\{(\gamma_i)\}_{i\in I}\subset \cone\big(\epi(f)-(x_0,f(x_0))\big)$ such that $\gamma_i\to -a=(0,-2),$ where
\begin{equation*}
\gamma_i=\mu_i\big((x_i,\lambda_i)-(x_0,f(x_0))\big),
\end{equation*}
where $\mu_i\geq 0,$ $(x_i,\lambda_i)\in\epi(f).$
It follows that
\begin{equation*}
\mu_i(x_i-x_0)\to 0\; \text{and}\; \mu_i(\lambda_i-f(x_0))\to-2.
\end{equation*}
Choose an index $i_0\in I$ such that
\begin{equation*}
  \mu_i(\lambda_i-f(x_0))<-1\; \text{whenever}\; i\geq i_0.
\end{equation*}
By the construction of $f$ we have the relationships
\begin{equation*}
\la x^*,\mu_i(x_i-x_0)\ra =\mu_i(f(x_i)-f(x_0)) \leq\mu_i(\lambda_i-f(x_0))<-1 \;\text{whenever}\;i\ge i_0
\end{equation*}
while contradicting the fact that $\la x^*,\mu_i(x_i-x_0)\to 0$. Thus $\lambda_0>f(x_0)$, which verifies \eqref{qrinotre} and completes the proof of the proposition.
$\h$

Finally, we consider epigraphs of extended-real-valued convex functions and present far-going LCTV space extensions of the finite-dimensional result of Corollary~\ref{epi-roc}.

\begin{Theorem} {\bf(quasi-relative interiors of convex epigraphs in LCTV spaces).}\label{quasi epi} Let $f\colon X\to\oR$ be a convex function, where $X$ is an LCTV space.  \begin{enumerate}
\item We have the inclusion
\begin{equation}\label{qriepi1}
{\rm qri}\big(\epi(f)\big)\supset\big\{(x,\lambda)\in X\times\R\;\big|\;x\in{\rm qri}\big(\dom(f)\big),\;\lambda>f(x)\big\}.
\end{equation}
\item If in addition  that $\epi(f)$ is quasi-regular, then $\dom(f)$ is quasi-regular and
\begin{equation}\label{qriepi}
{\rm qri}\big(\epi(f)\big)=\big\{(x,\lambda)\in X\times\R\;\big|\;x\in{\rm qri}\big(\dom(f)\big),\;\lambda>f(x)\big\}.
\end{equation}
\end{enumerate}
\end{Theorem}
{\bf Proof.} (a) Note that  \eqref{qriepi1} follows from \cite[Corollary~9(iii)]{ZO}.\\[1ex]
(b) Suppose that $\epi(f)$ is quasi-regular. Fix any $x\in \qri(\dom(f)$. Then we can choose $\lambda>f(x)$. By (a), we have $(x, \lambda)\in \qri(\epi f)=\iri(\epi(f))$. Now, we can apply Theorem \ref{LMregularepi} to see that $x\in \iri(\dom(f))$. Therefore, $\qri(\dom(f)\subset \iri(\dom(f))$, which verifies the quasi-regularity of $\dom(f)$. Let $F\colon X\tto \R$ be defined by $F(x):=[f(x), \infty)$. Then $\gph(F)=\epi(f)$, $\dom(F)=\dom(f)$, and for $x\in \dom(f)$, we have $\qri(F(x))=(f(x), \infty)$. Under the quasi-regularity of $\gph(F)=\epi(f)$, the inclusion ``$\subset$'' in \eqref{qriepi} follows directly from Theorem \ref{rocthm2}(a). $\h$\\

{\bf Acknowledgements.} The authors are indebted to Nicolas Hadjisavvas, Pedro P\'erez-Aros, Constantin Z\u{a}linescu, and two anonymous referees for their careful reading of the original version of the paper with helpful remarks and suggestions that allowed us to significantly improve the obtained results.

\end{document}